\documentclass[a4paper,11pt]{amsart}

\usepackage[T1]{fontenc}
\usepackage{amsthm}
\usepackage{amsmath}
\usepackage[authoryear]{natbib}
\usepackage[dvips,bookmarks=false]{hyperref}

\usepackage{setspace}
\allowdisplaybreaks[4]

\makeatletter

\RequirePackage{bbm}
\RequirePackage{mathrsfs}
\RequirePackage{amsfonts,amssymb}

\newcommand{\IN}{\mathbbm{N}}
\newcommand{\IZ}{\mathbbm{Z}}

\def\§#1{\mathscr{#1}}
\def\£#1{\mathcal{#1}}
\newcommand{\dist}{\§L}

\def\dtv{\mathop{d_{\mathrm{TV}}}}

\newcommand{\nsig}{\Sigma\kern-0.5em\raise0.2ex\hbox to 0pt{$\mid$}\kern0.5em}

\newcommand{\eps}{\varepsilon}
\newcommand{\D}{\Delta}

\newcommand{\ahalf}{{\textstyle\frac{1}{2}}}

\newcommand{\eq}{\eqref}
\newcommand{\IE}{\mathbbm{E}}
\newcommand{\IP}{\mathbbm{P}}
\newcommand{\Var}{\mathop{\mathrm{Var}}}

\newcommand{\e}{{\mathrm{e}}}
\newcommand{\Po}{\mathop{\mathrm{Po}}}

\def\be#1\ee{\begin{equation*}#1\end{equation*}}
\def\ben#1\ee{\begin{equation}#1\end{equation}}

\def\bs#1\es{\begin{split}#1\end{split}}
\def\bes#1\ee{\begin{equation*}\begin{split}#1\end{split}\end{equation*}}
\def\besn#1\ee{\begin{equation}\begin{split}#1\end{split}\end{equation}}

\def\bg#1\ee{\begin{gather*}#1\end{gather*}}
\def\bgn#1\ee{\begin{gather}#1\end{gather}}

\def\bm#1\ee{\begin{multline*}#1\end{multline*}}
\def\bmn#1\ee{\begin{multline}#1\end{multline}}

\def\ba#1\ee{\begin{align*}#1\end{align*}}
\def\ban#1\ee{\begin{align}#1\end{align}}

\def\bklr#1{\bigl(#1\bigr)}

\def\bkle#1{\bigl[#1\bigr]}

\def\klg#1{\{#1\}}
\def\bklg#1{\bigl\{#1\bigr\}}

\def\norm#1{\Vert#1\Vert}

\def\abs#1{\vert#1\vert}
\def\babs#1{\bigl\vert#1\bigr\vert}

\def\mid{\vert}
\def\bmid{\bigm\vert}

\def\^#1{\ifmmode {\mathaccent"705E #1} \else {\accent94 #1} \fi}
\def\~#1{\ifmmode {\mathaccent"707E #1} \else {\accent"7E #1} \fi}
\def\*#1{#1^\ast}
\def\>#1{\vec{#1}}
\def\.#1{\dot{#1}}

\def\leq{\leqslant}
\def\geq{\geqslant}

\numberwithin{equation}{section}
\theoremstyle{plain}
\newtheorem{theorem}{Theorem}[section]

\makeatother

\begin{document}

\title{A note on the exchangeability condition in Stein's method}
\author{Adrian R\"ollin}

\thanks{Supported by Swiss National Science projects 20--107935/1 and
PBZH2--117033} 

\maketitle

\begin{abstract}
\enskip We show by a surprisingly simple argument that the exchangeability
condition, which is key to the \emph{exchangeable pair} approach in Stein's
method for distributional approximation, can be omitted in many standard
settings. This is achieved by replacing the usual antisymmetric function by a
simpler one, for which only equality in distribution is required. In the case of
normal approximation we also slightly improve the constants appearing in
previous results. For Poisson approximation, a different antisymmetric function
is used, and additional error terms are needed if the bound is to be extended
beyond the exchangeable setting. 
\end{abstract}

\section{Introduction}

In the context of normal approximation, a variant of Stein's method that is
often used is the \emph{exchangeable pair coupling} introduced in
\cite{Diaconis1977} and \cite{Stein1986}. There are many applications based on
this coupling, see e.g.\ \cite{Rinott1997}, \cite{Fulman2004a},
\cite{Meckes2006a} and others, but also in the context of non-normal
approximation such as \cite{Chatterjee2005}, \cite{Chatterjee2006} and
\cite{Rollin2007a}.

Assume that $W$ is a random variable which we want to approximate. The key
concept introduced by \cite{Stein1986} is that of \emph{auxiliary
randomization}. In the context of exchangeable pairs, this means that we
construct a random variable $W'$ on the same probability space such that
$(W,W')$ is exchangeable, that is $\dist(W,W') = \dist(W',W)$ and, in general,
the pair should be constructed such that $\abs{W'-W}$ is small. One can then
prove for instance a bound as in Theorem~\ref{th1} below (often under additional
conditions on the exchangeable pair such as Equation~\eq{5}). If $W$ and $W'$
are two consecutive steps of a reversible Markov chain in equilibrium,
\cite{Rinott1997} note that exchangeability automatically follows. However, the
pairs constructed in some of the examples of \cite{Rinott1997} and
\cite{Fulman2004} are based on non-reversible Markov chains and therefore some
effort has to be put into showing that the pairs satisfy the exchangeability
condition. 

The key fact to prove a bound as in Theorem \ref{th1} is that, for any
antisymmetric function $F$, exchangeability of $(W,W')$ implies the identity 
\ben								\label{1}
 	\IE F(W,W') = 0
\ee
(see \cite[p.\,10]{Stein1986}). In fact, this is often the only place where
exchangeability is used. For example, in the case of the normal distribution,
the standard choice for $F$ is  
\ben								\label{2}
	F(w,w') = (w'-w)\bklr{f(w')+f(w)},
\ee
where $f$ is the solution to the Stein equation
\ben								\label{3}
	f'(x) - x f(x) = h(x) - \IE h(Z)
\ee
and $Z$ has the standard normal distribution. We show in this paper that instead
of the choice \eq{2} we can take the simpler function  
\ben								\label{4}
	F(w,w') = \int_0^{w'} f(x) d x - \int_0^{w} f(x) d x,	
\ee
which is, in particular, antisymmetric, but for which \eq{1} of course holds
without exchangeability as long as $\dist(W') = \dist(W)$. If exchangeable pairs
are used in a discrete setting, the integrals in \eq{4} of course have to be
replaced by corresponding sums.

In the following section we restate and prove results obtained by the
exchangeable pairs approach for normal and Poisson approximation. The
main purpose is to demonstrate in detail how the exchangeability condition can
be omitted in some standard settings, sometimes also yielding better constants.
With the approach in this paper, however, it is possible to remove the
exchangeability condition in many other situations, such as in \cite{Stein1995}
and \cite{Meckes2006a}, \cite{Chatterjee2006}, and \cite{Rollin2007a}.

Most ingredients in the following proofs are taken directly from the proofs of
the corresponding papers; hence many details have been omitted.  

\section{Main results}

\subsection{Normal approximation}

For normal approximation we need some more assumptions. Assume that $W$ is a
random variable with $\IE W = 0$ and $\Var W=1$. Construct an exchangeable pair
such that 
\ben								\label{5}
	\IE^W W' = (1-\lambda)W + R
\ee 
for some constant $0<\lambda<1$ and some random variable $R$, where $\IE^W$
denotes the conditional expectation with respect to $W$. 

The following calculations are essential for the proof of
Theorem~1.2 of \cite{Rinott1997}, but we now assume for the sake of clarity
that $R=0$. Let $f$ be the solution to \eq{3} for a Lipschitz-continuous test
function~$h$. Now, with the standard choice of $F$ as
in \eq{2} and exchangeability, we have 
\besn								\label{6}
	0 & =  \IE F(W,W') \\
	& = \IE\bklg{(W'-W)\bklr{f(W')+f(W)}}\\
	& = \IE\bklg{(W'-W)\bklr{2f(W) + f(W')-f(W)}}\\
	& = -2\lambda\IE\bklg{ W f(W)} + \IE\bklg{(W'-W)\bklr{f(W')-f(W)}},
\ee
where for the last equality we used \eq{5}. Let now $\tau$ be a random variable
uniformly distributed on $[0,1]$,  independent of all other random variables,
and put $V = W'-W$. Noting that the second derivative $f''$ exists almost
everywhere because $h$  is Lipschitz continuous, Taylor's expansion yields 
\be
	f(W') = f(W) + V f'(W) + V^2\IE^{W,W'}\bklg{(1-\tau)f''(W+\tau V)}.
\ee
Thus, from \eq{6},
\ben								\label{7}
    \IE\bklg{W f(W)} = \frac{1}{2\lambda}\IE\bklg{V^2 f'(W)} 
    + \frac{1}{2\lambda}\IE\bklg{V^3(1-\tau)f''(W+\tau V)}.
\ee

In the proof of Theorem \ref{th1} we show that an equation similar to \eq{7} can
be deduced without exchangeability. To state the theorem, we need some notation.
Define for a given function $h$ and $\eps>0$,
\be
	h^+_\eps(x)= \sup\{h(x+y)\,:\,\abs{y}\leq\eps\},
    \quad h^-_\eps(x)= \inf\{h(x+y)\,:\,\abs{y}\leq\eps\}.
\ee
Let $\£H$ be a class of measurable functions on the real line such that 
for all $h\in\£H$, we have $\norm{h}\leq 1$, where $\norm{\!\cdot\!}$
denotes the supremum norm; 
for any real numbers $c$ and $d$, $h\in\£H$ implies $h(c\cdot +d)\in\£H$; 
for any $\eps>0$, $h\in\£H$ implies $h^+_\eps,h^-_\eps\in\£H$ and there is
a constant $a$ (depending only on the class $\£H$) such that 
$\IE\bklg{h^+_\eps(Z) - h^-_\eps(Z)}\leq a\eps$ where $Z$ has standard normal
distribution.
As in \cite{Rinott1997} we assume without loss of generality that
$a\geq\sqrt{2/\pi}$. 

\begin{theorem}[cf. Theorem 1.2 of \cite{Rinott1997}]\label{th1} Assume that $W$
and $W'$ are random variables on the same probability space such that
$\dist(W')=\dist(W)$, $\IE W = 0$, $\Var W = 1$. Assume that \eq{5} holds for
some $\lambda$ and $R$. Then, for $\delta := \sup_{h\in\£H}\babs{\IE
h(W)-\IE h(Z)}$, 
\ben								\label{8}
	\delta \leq \frac{6}{\lambda}\sqrt{\Var \IE^W(W'-W)^2}
    +\frac{19\sqrt{\IE R^2}}{\lambda} + 
		4\sqrt\frac{a\IE\abs{W'-W}^3}{\lambda}.
\ee
If, in addition, there is a constant $A$ such that $\abs{W'-W}\leq A$ almost
surely, we have 
\ben								.\label{9}
	\delta \leq \frac{12}{\lambda}\sqrt{\Var \IE^W (W'-W)^2}
	+ \frac{37\sqrt{\IE R^2}}{\lambda}
	+ 32\frac{A^3}{\lambda} 
	+ 6\frac{A^2}{\sqrt{\lambda}}
\ee
\end{theorem}

\begin{proof} From Lemma 4.1 of \cite{Rinott1997} we have that, for any $0<t<1$,
\ben								\label{10}
	\delta \leq 2.8\sup_{h\in\£H}\babs{\IE h_t(W)-\IE h_t(Z)} + 4.7 a t,
\ee
where $h_t(x) = \IE h(x + t Z)$. Let $f$ be the solution to the Stein equation
\ben                                        \label{11}
	f'(x) - x f(x) = h_t(x) - \IE h_t(Z).
\ee
Then, $f$ satisfies
\ben								\label{12}
	\norm{f}\leq 2.6,\quad\norm{f'}\leq 4,
    \quad\norm{f''}\leq2\norm{h'_t}\leq 1.6t^{-1};
\ee
recalling that $\norm{h}\leq 1$, the first two bounds and the first part of
the third one follow from Lemma~3 of \cite{Stein1986} and, as noted by
\cite{Rinott1997}, the second inequality of the third bound can be deduced using
the equality $h'_t(x)=-t^{-1}\int h(x+t y)\varphi'(y) d y$, where $\varphi$ is
the standard normal density, so that $\norm{h'_t}\leq
t^{-1}\int\abs{\varphi'(x)}dx = t^{-1}\sqrt{2/\pi}$.
 
Define the function 
\ben								\label{13}
	G(w) = \int_0^w f(x)\,d x
\ee
and note that $\abs{G(w)}\leq \abs{w}\norm{f}$, so that $\IE G(W)$ exists. By
Taylor's expansion, we have 
\besn								\label{14}
	G(W') & = G(W) + V f(W) + \ahalf V^2 f'(W) \\
	& \qquad + \ahalf V^3\IE^{W,W'}\bklg{(1-\tau)^2 f''(W+\tau V)},
\ee
where, again, $V=W'-W$. Thus, together with \eq{5}, we obtain
\bes
	0& = \IE G(W') - \IE G(W) \\
	& = -\lambda\IE\bklg{W f(W)}+ \IE\bklg{R f(W)}  \\
	& \quad+ \ahalf\IE\bklg{V^2 f'(W)} 
		+ \ahalf\IE\bklg{ V^3(1-\tau)^2 f''(W+\tau V)},
\ee
which can be rearranged to obtain the following analogue of \eq{7}
\besn								\label{15}
	\lambda\IE\bklg{W f(W)} & 
    =\ahalf \IE\bklg{V^2 f'(W)}
    +\ahalf\IE\bklg{V^3(1-\tau)^2f''(W+\tau V)} \\ 
	&\quad + \IE\bklg{R f(W)}.
\ee
With $\alpha:=\IE\klg{R W}$ and noting that $\IE V^2 = 2(\lambda-\alpha)$ from
\eq{5}, we thus have from \eq{11} and \eq{15} 
\ba
	\lambda\bklr{\IE h_t(W) - \IE h_t(Z)}
	& = \IE \bklg{\bklr{(\lambda-\alpha)-\ahalf V^2}f'(W)}\\
	&\quad + \IE\bklg{\alpha f'(W)-R f(W)}\\
	&\quad - \ahalf\IE\bklg{V^3(1-\tau)^2f''(W+\tau V)}\\
	& =: J_1 + J_2 - \ahalf J_3.
\ee
Using \eq{12}, we obtain the estimates
\be
	\abs{J_1} \leq 2\sqrt{\Var \IE^W V^2},
    \qquad \abs{J_2} \leq 6.6\sqrt{\IE R^2},
    \qquad\abs{J_3} \leq \frac{1.6}{3t}\IE\abs{V}^3;
\ee
for details see \cite{Rinott1997}. Choosing
$t=0.4\bklr{\IE\abs{V}^3/(a\lambda)}^{1/2}$ and with \eq{10}, this proves
\eq{5}. 

Assume now that $\abs{V}\leq A$. Note that because of \eq{11}, $f''(x) = f(x) +
x f'(x) + h'_t(x)$. Following the proof of \cite{Rinott1997}, but recalling that
in our remainder $J_3$ the term $(1-\tau)$ is squared, we have  
\bes
	\abs{J_3} & = \IE\bklg{V^3(1-\tau)^2\bklr{f(W+\tau V) 
        + (W+\tau V)f'(W+\tau V) + h'_t(W+\tau V)}}\\
	& \leq 0.9 A^3 + 1.4 \IE\bklg{V^3(\abs{W}+\abs{W'})} 
        + \IE\bklg{V^3(1-\tau)^2h'_t(W+\tau V)}\\
	& \leq 3.7 A^3 + \IE\bklg{V^3(1-\tau)^2h'_t(W+\tau V)},
\ee
where the latter expectation can be bounded by
\be
	\babs{\IE\bklg{V^3(1-\tau)^2h'_t(W+\tau V)} }
    \leq \frac{a A^3}{3} + \frac{A^3}{3t}(2\delta + a A) ;
\ee
see again \cite{Rinott1997}. Collecting all the bounds on the $J_i$ we obtain
\besn								\label{16}
	\babs{\IE h_t(W)-\IE h_t(Z)} 
		& \leq \frac{2}{\lambda}\sqrt{\Var\IE^WV^2}
            + \frac{6.6}{\lambda}\sqrt{\IE R^2}\\
		&\quad + \frac{3.7A^3}{2\lambda} + \frac{a A^3}{6\lambda} 
            + \frac{A^3(2\delta+a A)}{6\lambda t}.
\ee
Recalling that $a\geq\sqrt{2/\pi}$, putting \eq{16} into \eq{10} and with the
choice
$t=0.32 A \bklr{\frac{A(2\delta+a A)}{a\lambda}}^{1/2}$, we finally have
\be
	\delta 
	\leq 5.6\sqrt{\Var \IE^W V^2} 
	+ 18.5 \sqrt{\IE R^2} + 7 \frac{a A^3}{\lambda}
	+ 3 \frac{a A^2}{\sqrt{\lambda}} 
	+ 4.2\sqrt{\frac{a\delta A^3}{\lambda}}.
\ee
This inequality is of the form $\delta\leq a+b\sqrt{\delta}$, for which we can
show that $\delta\leq 2a+b^2$ and which hence proves \eq{9}. 
\end{proof}

Note that if $R=0$ almost surely, equation \eq{15} is the same as \eq{7}, except
that the factor $(1-\tau)$ is squared; this is the only reason for the improved
constants. If $(W,W')$ is exchangeable, both equalities \eq{7} and \eq{15} hold.
At first glance this may seem to be a contradiction, but the remainders with the
second derivatives are in fact equal. To see this, write 
\bes
	&\IE\bklg{V^3(1-\tau)^2 f''(W+\tau V)}\\
	&\qquad = \IE\bklg{V^3(1-\tau) f''(W+\tau V)} 
        - \IE\bklg{V^3\tau(1-\tau) f''(W+\tau V)}.
\ee
We need only show that the second term on the right hand side is equal to zero.
Note to this end the simple fact that $\dist(\tau) = \dist(1-\tau)$, thus, using
this in the following calculations to obtain the first equality and the
exchangeability of $(W,W')$ for the third equality, 
\ba
	&\IE\bklg{(W'-W)^3\tau(1-\tau) f''\bklr{W+\tau (W'-W)}} \\
	&\qquad=\IE\bklg{(W'-W)^3\tau(1-\tau) f''\bklr{W+(1-\tau) (W'-W)}}\\
	&\qquad=\IE\bklg{(W'-W)^3\tau(1-\tau) f''\bklr{W'+\tau (W-W')}}\\
	&\qquad=\IE\bklg{(W-W')^3\tau(1-\tau) f''\bklr{W+\tau (W'-W)}}\\
	&\qquad=-\IE\bklg{(W'-W)^3\tau(1-\tau) f''\bklr{W+\tau (W'-W)}},\\
\ee
which proves the claim.

\subsection{Poisson approximation}

The situation in the discrete setting, in which $W$ takes values only in $\IZ$,
is more delicate. Assume that $f$ is the solution to a Stein equation of a
discrete distribution. Instead of \eq{13}, define the function $G$ as 
\ben								\label{17}
	G(w) = \sum_{k=1}^w f(k) - \sum_{k=0}^{-w-1} f(-k),
\ee
where here and in what follows $\sum_{k=a}^b$ is defined to be zero if $b<a$.
For many standard distributions, $f$ will be bounded, thus $\abs{G(w)}\leq
\abs{w}\norm{f}$, so that $\IE G(W)$ exists if $\IE\abs{W}<\infty$. 

With such $G$, one verifies that $G(w)-G(w-1) = f(w)$ for all $w\in\IZ$. 
Define $I_i = I[W'-W = i]$ and $\D_i G(w) := G(w+i)-G(w)$ for all $i\in\IZ$;
then, 
\besn								\label{18}
	G(W') - G(W) & = \sum_{i\in\IZ} I_i \Delta_i G(W)
\ee
and thus, if $\dist(W') = \dist(W)$, 
\besn								\label{19}
	0&=\IE G(W') - \IE G(W) 
    = \sum_{i\in\IZ} \IE\bklg{I_i\Delta_i G(W)}
    =\sum_{i\in\IZ}\IE\bklg{ P_i(W) \Delta_i G(W)}
\ee
where $P_i(W) := \IP^W[W'-W = i]$.

As mentioned in the introduction, \cite{Rinott1997} suggest constructing $W$ and
$W'$ via an underlying stationary Markov process, which implies
$\dist(W)=\dist(W')$; if the chain is reversible, exchangeability follows
immediately. However, the Markov chain that they use for the anti-voter model is
not reversible, and so exchangeability has to be proved separately. One way of
doing this is to assume that  
\ben								\label{20}
    W'-W\in\{-1,0,1\},
\ee
almost surely, which is the main assumption in \cite{Rollin2007a}. 

\cite{Chatterjee2005} applied the method of exchangeable pairs to the Poisson
distribution. Assume that $(W',W)$ is an exchangeable pair taking values only on
the non-negative integers. Define the antisymmetric function 
\ben								\label{21}
	F(w,w') = f(W') I[W'-W = 1] - f(W) I[W'-W = -1].
\ee
This yields
\besn								\label{22}
	0 & = \IE F(W,W')\\
	& = \IE\bklg{f(W') I[W'-W = 1] - f(W) I[W'-W = -1]} \\
	& = \IE\bklg{\IE^W\klg{f(W') I[W'-W = 1]} 
		- \IE^W\klg{f(W) I[W'-W = -1]}}\\
	& = \IE\bklg{f(W+1)\IP[W'-W = 1|W] - f(W)\IP[W'-W = -1|W] }\\
	& =\IE\bklg{f(W+1)P_{1}(W) - f(W)P_{-1}(W) }.
\ee 
We can use now the standard argument for Poisson approximation by Stein's method
(see \cite[p.~6]{Barbour1992}). Denoting by $\Po(\lambda)$ the Poisson
distribution with mean $\lambda$, it follows that
\besn								\label{25}
	&\dtv\bklr{\dist(W),\Po(\lambda)} 
    \leq \sup_f \babs{\IE\bklg{\lambda f(W+1) - W f(W)}}\\ 
	&\leq \sup_f\babs{\IE\bklg{(c P_{1}(W)-\lambda)f(W+1) 
    - (c P_{-1}(W)-W)f(W)}} =: \kappa_c
\ee
for any $c>0$ and where the supremum ranges over all solutions $f=f_A$ to the
Stein-equation
\ben								\label{26}
	\lambda f(j+1) - j f(j) = I[j\in A] - \Po(\lambda)\klg{A},
\ee
for subsets $A$ of the non-negative integers. For many applications investigated
by \cite{Chatterjee2005}, the following further bound on $\kappa_c$ is used:
\be
	\kappa_c 
	\leq \lambda^{-1/2}\bklr{\IE\abs{c P_1(W)-\lambda}
		+\IE\abs{cP_{-1}(W)-W}},
\ee
which results from the well known bound $\norm{f}\leq \lambda^{-1/2}$ (see
\cite[Lemma 1.1.1]{Barbour1992}). However, it is at times beneficial to work
directly with the expression in \eq{25} in the hope of re-arranging things so
that the better bound $\norm{\D_1 f}\leq (1-\e^{-\lambda})\lambda^{-1}$ may be
applied.

It may seem surprising that, although in the bound \eq{25} only the jump
probabilities of size $1$ appear, no assumptions concerning the sizes of other
jumps were made in the above calculations; we did not, for instance, assume
condition \eq{20}. However, with the choice $f(w) = I[w = k]$ for $k\in\IZ$, we
obtain from \eq{21} the detailed balance equation for reversible Markov chains.
Even if the chain makes jumps of size larger than $1$, the stationary
distribution is determined by the jump probabilities of size $1$, so that in
\eq{25} the full information about the distribution under consideration is
actually used, just by starting from \eq{21}. If exchangeability is not assumed,
the effects of jumps of size larger that $1$ have also to enter, and this is
reflected in the bounds of the next theorem. 
The choice of $F$ in \eq{21} is fundamentally different from the standard choice
in the continuous setting, where we obtained the same bounds as before, but
under weaker assumptions.

\begin{theorem}[cf.\ Proposition 3 of \cite{Chatterjee2005}]\label{th2} Let $W$
and $W'$ be non-negative random variables such that $\dist(W')=\dist(W)$. Then,
for any constant $c>0$, 
\ben								\label{27}	
	\dtv\bklr{\dist(W),\Po(\lambda)} \leq \kappa_c + c\rho
\ee
where $\rho$ satisfies the bounds
\ben								\label{28}
	\rho  
	\leq\lambda^{-1/2}
		\sum_{i\geq 2}i\sum_{k\in\IZ}\babs{p_{k,k+i}-p_{k+i,k}} 
	\leq \lambda^{-1/2}\sum_{\abs{i}\geq 2}\abs{i}\IE P_i(W),
\ee
for $p_{k,j} = \IP[W=k,W'=j]$ and the $P_i$ are as before.

\end{theorem}
\begin{proof}
Taking expectation over \eq{26} with respect to $W$, using \eq{19} and noting
that $\D_1 G(W) = f(W+1)$ and $\D_{-1} G(W) = - f(W)$, we obtain
\bes
	&\IE\bklg{\lambda f(W+1)-W f(W)}\\
	&\qquad = \IE\bklg{\bklr{\lambda- c P_1(W)} f(W+1)
			- \bklr{W - c P_{-1}(W)}f(W)}\\
	&\qquad\quad+c\sum_{i\geq 2} \IE\bklg{I_i\D_i G(W) 
		+ I_{-i}\D_{-i} G(W)},
\ee
where, as before, $I_i = I[W'-W=i]$. Now, it is easy to see that
\ben								\label{29}
	\IE\bklg{I_i\D_i G(W) + I_{-i}\D_{-i} G(W)} 
	= \sum_{k\in\IZ}(p_{k,k+i}-p_{k+i,k})\D_i G(k)
\ee
Recalling the bound $\norm{f}\leq \lambda^{-1/2}$ for all solutions $f$ of
\eq{26}, and hence $\abs{\Delta_i G(W)}\leq\abs{i}\lambda^{-1/2}$ for all
$i\in\IZ$, \eq{29} yields the first bound of \eq{28}. Now, as $p_{k,k+i} =
\IP[W=k]P_i(k)$, the second bound follows from the first.
\end{proof}

Under condition \eq{20}, Theorem \ref{th2} and estimate \eq{25} clearly yield
the same bound. If \eq{20} does not hold, exchangeability is not automatically
implied. Then the first bound of \eq{28} is a measure of the non-exchangeability
of $(W,W')$; if the pair is in fact exchangeable, this term vanishes and we
regain \eq{20}. The second bound in \eq{28} is particularly useful if the jump
probabilities of larger jumps are small.

\section{Discussion}

The key to understand the present approach is the generator interpretation
introduced in \cite{Barbour1988}. Recall that
\be
	(\£A G)(x) = G''(x) - x G'(x)	= f'(x) - x f(x)
\ee
is the generator of the Ornstein-Uhlenbeck diffusion, applied to the
function~$G$. Now assume that a Markov chain $\{X(n); n\in\IZ_+\}$ with
stationary distribution $\dist(W)$ is given and let $(W,W')=(X(0),X(1))$,
assuming that the chain starts in its equilibrium. Construct a Markov jump
process $Z = \{Z(t); t\geq 0\}$ by randomising the fixed steps of size $1$ of
the Markov chain $X$ and wait instead an exponentially distributed amount of
time with rate $1/\lambda$. It is easy to see that, under the assumption of
\eq{5}, the infinitesimal operator of $Z$ is
\ba
	&(\£B G)(x)   
	= \lim_{h\to 0}\frac{\IE\bkle{G(Z(h))\bmid Z(0)=x} - G(x)}{h}
	= \frac{\IE^x G(W') - G(x)}{\lambda}\\
	&\quad = -x G'(x) + \frac{\IE^x V^2}{2\lambda} G''(x) \\
	&\qquad\qquad\qquad	+
	\int\frac{\IE^x\bklg{(V-t)^2 I[\text{$V>t>0$ or  $V<t<0$}]}}{2\lambda}
			\,G^{(3)}(x+t) d t 
\ee
where $\IE^x$ denotes $\IE\bklg{\cdot\,|W=x}$ and $V=W'-W$. Thus, we in fact
compare the infinitesimal operator of this jump process with the generator of
the Ornstein-Uhlenbeck diffusion. 

The antisymmetric function $F(w,w') = (w'-w)(f(w')+f(w))$ in fact also (almost)
calculates this infinitesimal operator, but with the first step of the Taylor's
expansion already carried out.

In the case of \cite{Stein1995} and \cite{Meckes2006a}, we have a family of
Markov jump processes $\{Z_m;m\in\IN\}$ with infinitesimal operators
\be
	(\£B_m G)(x) = \frac{ \IE^x G(W_m) - G(x)}{\eps_m\vartheta}
\ee 
where $\eps_m^{-1}\lambda^{-1}$ are the jump rates and by letting $\eps\to0$ we
show that the processes $Z_m$ converge to a diffusion process with infinitesimal
operator 
\be	
	(\£B G)(x) = (1+\ahalf\lambda^{-1}\IE^x E) G''(x)-x G'(x) 
\ee
for some specific random variable $E$. This diffusion has the same linear drift
as the Ornstein-Uhlenbeck diffusion, but a non-constant diffusion rate. Stein's
method now allows us to state through the approaches in \cite{Stein1995} and
\cite{Meckes2006a} that the less $\IE^x E$ fluctuates around zero, the nearer
the stationary distribution of the process is to the stationary distribution of
Ornstein-Uhlenbeck diffusion, that is, to the standard normal distribution.

\section*{Acknowledgements}

I thank Andrew Barbour for helpful discussions and the anonymous referee for
useful comments on the manuscript.


\begin{thebibliography}{12}
\providecommand{\natexlab}[1]{#1}
\providecommand{\url}[1]{\texttt{#1}}
\expandafter\ifx\csname urlstyle\endcsname\relax
  \providecommand{\doi}[1]{doi: #1}\else
  \providecommand{\doi}{doi: \begingroup \urlstyle{rm}\Url}\fi

\bibitem[Barbour(1988)]{Barbour1988}
A.~D. Barbour (1988).
\newblock Stein's method and {P}oisson process convergence.
\newblock \emph{J.~Appl. Probab.} \textbf{25A}, \penalty0 175--184.

\bibitem[Barbour et~al.(1992)Barbour, Holst, and Janson]{Barbour1992}
A.~D. Barbour, L.~Holst, and S.~Janson (1992).
\newblock \emph{Poisson approximation}, volume~2 of \emph{Oxford Studies in
  Probability}.
\newblock The Clarendon Press Oxford University Press, New York.
\newblock ISBN 0-19-852235-5.
\newblock Oxford Science Publications.

\bibitem[Chatterjee and Fulman(2006)]{Chatterjee2006}
S.~Chatterjee and J.~Fulman (2006).
\newblock Exponential approximation by exchangeable pairs and spectral graph
  theory.
\newblock \emph{Preprint}.
\newblock Available at \url{http://arxiv.org/abs/math/0605552}.

\bibitem[Chatterjee et~al.(2005)Chatterjee, Diaconis, and
  Meckes]{Chatterjee2005}
S.~Chatterjee, P.~Diaconis, and E.~Meckes (2005).
\newblock Exchangeable pairs and {P}oisson approximation.
\newblock \emph{Probab. Surv.} \textbf{2}, \penalty0 64--106.
\newblock Available at \url{www.arxiv.org/abs/math.PR/0411525}.

\bibitem[Diaconis(1977)]{Diaconis1977}
P.~Diaconis (1977).
\newblock The distribution of leading digits and uniform distribution {${\rm
  mod}$}~{$1$}.
\newblock \emph{Ann. Probability} \textbf{5}, \penalty0 72--81.

\bibitem[Fulman(2004{\natexlab{a}})]{Fulman2004}
J.~Fulman (2004{\natexlab{a}}).
\newblock Stein's method and non-reversible {M}arkov chains.
\newblock In \emph{Stein's method: expository lectures and applications},
  volume~46 of \emph{IMS Lecture Notes Monogr. Ser.}, pages 69--77. Inst. Math.
  Statist.

\bibitem[Fulman(2004{\natexlab{b}})]{Fulman2004a}
J.~Fulman (2004{\natexlab{b}}).
\newblock Stein's method, {J}ack measure, and the {M}etropolis algorithm.
\newblock \emph{J. Combin. Theory Ser. A} \textbf{108}, \penalty0 275--296.

\bibitem[Meckes(2006)]{Meckes2006a}
E.~Meckes (2006).
\newblock Linear functions on the classical matrix groups.
\newblock \emph{Preprint}.
\newblock Available at \url{www.arxiv.org/math.PR/0509441}.

\bibitem[Rinott and Rotar(1997)]{Rinott1997}
Y.~Rinott and V.~Rotar (1997).
\newblock On coupling constructions and rates in the {CLT} for dependent
  summands with applications to the antivoter model and weighted
  {$U$}-statistics.
\newblock \emph{Ann. Appl. Probab.} \textbf{7}, \penalty0 1080--1105.

\bibitem[R\"ollin(2007)]{Rollin2007a}
A.~R\"ollin (2007).
\newblock Translated {P}oisson approximation using exchangeable pair couplings.
\newblock \emph{Ann. Appl. Prob.} \textbf{17}, \penalty0 1596--1614.

\bibitem[Stein(1986)]{Stein1986}
C.~Stein (1986).
\newblock \emph{Approximate computation of expectations}.
\newblock Institute of Mathematical Statistics Lecture Notes---Monograph
  Series, 7. Institute of Mathematical Statistics, Hayward, CA.

\bibitem[Stein(1995)]{Stein1995}
C.~Stein (1995).
\newblock The accuracy of the normal approximation to the distribution of the
  traces of powers of random orthogonal matrices.
\newblock Technical Report No.~470, Standford University Department of
  Statistics.

\end{thebibliography}
\end{document}